\documentclass[fleqn,12pt]{article}

\usepackage[cp1251]{inputenc}
\usepackage{srcltx}
\usepackage{amsfonts}
\usepackage{amsmath}

\newcommand{\rank}{{\rm rank}}
\begin{document}

\fontsize{14.4}{18pt}\selectfont

\noindent UDC 517.9 \vspace{6.0mm}

\noindent \textbf{A. M. Samoilenko}\quad (Institute of Mathematics, Ukrainian
National Academy of Sciences, Kyiv)\vspace{6.0mm}

\noindent \textbf{ON INVARIANT MANIFOLDS OF LINEAR DIFFERENTIAL EQUATIONS}
\vspace{6.0mm}

\noindent In the present paper, we develop the ideas of the Bogolyubov method
of integral manifolds for linear differential equations. The obtained result
can also have a practical interpretation and be used in applications.

\vspace{6.0mm}

\begin{center}\textbf{ 1. Two Lemmas}\end{center}

\rm Suppose that matrices  $\Phi_{1}(t)$ and $\Phi_{2}(t)$ are continuously
differentiable for all   $t\in \mathbb{R}$, $\Phi_{1}(t)\in
M_{nm}(\mathbb{R})$, $\Phi_{2}(t)\in M_{pm}(\mathbb{R}),$ where $m>n>p=m-n,$
and the following condition is satisfied:
\begin{equation}\det\left(\begin{array}{cc}\Phi_{1}(t)\\\Phi_{2}(t)\end{array}\right)\ne 0.
\nonumber \end{equation}

Let $\Phi_{1}^{+}(t)$ and $\Phi_{2}^{+}(t)$ denote matrices pseudoinverse to
$\Phi_{1}(t)$ and $\Phi_{2}(t)$ and defined by the condition
\begin{equation}
\label{1}
\left(\begin{array}{cc}\Phi_{1}(t)\\\Phi_{2}(t)\end{array}\right)(\Phi_{1}^{+}(t),\Phi_{2}^{+}(t))=E,
\,\,\,\,E\in M_{m}(\mathbb{R}).
\end{equation}
Let
\begin{equation}   \label{2}
M_{1}(t)=\Phi_{1}^{+}(t)\Phi_{1}(t),
\end{equation}
\begin{equation}   \label{3}
M_{2}(t)=\Phi_{2}^{+}(t)\Phi_{2}(t).
\end{equation}
If follows from (\ref{1}) and the definitions of the matrices $M_{1}(t)$ and
$M_{2}(t)$ that these matrices satisfy the conditions
\begin{equation}
M_{\nu}^{2}(t)=M_{1}(t),\,\,\,\,\nu=1,2, \qquad \rank\,M_{1}(t)=n,
\qquad \rank\,M_{2}(t)=p,\nonumber
\end{equation}
\begin{equation}   \label{4}
M_{1}(t)M_{2}(t)=M_{2}(t)M_{1}(t)=0, \qquad M_{1}(t)+M_{2}(t)=E.
\end{equation}
Let us prove equality (\ref{4}). Indeed, since the matrix
$(\Phi_{1}^{+}(t),\Phi_{2}^{+}(t))$ is inverse to the matrix
$\left(\begin{array}{cc}\Phi_{1}(t)\\\Phi_{2}(t)\end{array}\right),$
multiplying the latter from the left by the matrix
$(\Phi_{1}^{+}(t),\Phi_{2}^{+}(t))$  we obtain the equality
$$\
Phi_{1}^{+}(t)\Phi_{1}(t)+\Phi_{2}^{+}(t)\Phi_{2}(t)=E,
$$
which is equivalent to equality (\ref{4}).

In the space $\mathbb{R}^{m},$ we define two subspaces by using the matrix
$M_{1}(t)$, namely,
\begin{equation}   \label{5}
M^{n}(t)=\{y\in \mathbb{R}^{m}: y=M_{1}(t)y\},
\end{equation}
\begin{equation}   \label{6}
M^{m-n}(t)=\{y\in \mathbb{R}^{m}: M_{1}(t)y=0\},
\end{equation}
and two subspaces by using the matrix $M_{2}(t)$:
\begin{equation}   \label{7}
M_{1}^{p}(t)=\{y\in \mathbb{R}^{m}: y=M_{2}(t)y\},
\end{equation}
\begin{equation}   \label{8}
M_{1}^{m-p}(t)=\{y\in \mathbb{R}^{m}: M_{2}(t)y=0\}.
\end{equation}

\textbf{Lemma 1.}\,\,\,  {\em The subspaces  \/} $M^{k}(t)$ {\em and \/}
$M_{1}^{k}(t),\,\,\,\,k\in\{n,p\}$, {\em satisfy the conditions \/}
\begin{equation} M^{n}(t)=M_{1}^{m-p}(t)=\ker \,\Phi_{2}(t),\,\,\,\,M^{m-n}(t)=M_{1}^{p}(t)=
\ker \,\Phi_{1}(t).  \nonumber
\end{equation}

Indeed, according to properties of the matrix $M_{1}(t),$ the general solution of the
equation defined by subspace (\ref{5}) is the function
\begin{equation}   \label{9}
y=M_{1}(t)c(t),
\end{equation}
where $c(t)$ is an arbitrary function with values in $\mathbb{R}^{m}.$ According to
properties of the matrix  $M_{2}(t),$ the general solution of the equation
defined by subspace (\ref{8}) is the function
\begin{equation}   \label{10}
y=(E-M_{2}(t))c_{1}(t)=M_{1}(t)c_{1}(t).
\end{equation}
Thus, the general solutions of the considered equations coincide for
$c(t)=c_{1}(t),$ which proves the equality $M^{n}(t)=M_{1}^{m-p}(t).$

It follows from the definition of  $M_{1}^{m-p}$ that
\begin{equation}
M_{1}^{m-p}(t)=\ker \,M_{2}(t)=\ker \,(\Phi_{2}^{+}(t)\Phi_{2}(t)). \nonumber
\end{equation}
Taking into account that the equality $\Phi_{2}(t)x=0$ implies that  $\Phi_{2}(t)
\Phi_{2}^{+}(t)x=x=0,$ we conclude that $\ker \,\Phi_{2}^{+}(t)=0$ and, hence, $\ker
\,(\Phi_{2}^{+}(t)\Phi_{2}(t))=\ker \,\Phi_{2}(t).$ This proves the first
equality of Lemma~1. The second equality is proved by analogy.

\textbf{Lemma 2.}\,\,\, {\em The mapping \/} $\Phi_{1}^{+}(t):y=\Phi_{1}^{+}(t)x$
{\em is a diffeomorphism  of \/} $\mathbb{R}^{n}$ {\em into \/} $M^{n}(t)$, {\em
and the mapping \/} $\Phi_{2}^{+}(t):y=\Phi_{2}^{+}(t)x$ {\em is a diffeomorphism
of \/} $\mathbb{R}^{n}$ {\em into \/} $M^{m-n}(t).$

We prove only the first assertion of Lemma~2 because the second assertion is proved by analogy.

The matrix $\Phi_{1}^{+}(t),$ as a block of the matrix inverse to
$\left(\begin{array}{cc}\Phi_{1}(t)\\\Phi_{2}(t)\end{array}\right),$ is
continuously differentiable and has a continuously differentiable pseudoinverse
matrix, namely  $\Phi_{1}(t).$ Therefore, to prove the first assertion of
Lemma~2, it remains to prove that the mapping
$\Phi_{1}^{+}(t):\mathbb{R}^{n}\rightarrow M^{n}(t)$ is a homeomorphism.

Since $\ker \Phi_{1}^{+}(t)=0,$ we conclude that $\Phi_{1}^{+}(t)$ is a homeomorphism of
$\mathbb{R}^{n}$ into the image of $\Phi_{1}^{+}(t;\mathbb{R}^{n})$ under the
mapping $\Phi_{1}^{+}(t):\mathbb{R}^{n}\rightarrow \mathbb{R}^{m}.$ It remains
to prove that $\Phi_{1}^{+}(t;\mathbb{R}^{n})=M^{n}(t).$

Assume that this is not true. Then either there exists a point
$y\in\Phi_{1}^{+}(t;\mathbb{R}^{n})$ such that  $y\notin M^{n}(t)$ or there
exists a point  $y\in \mathbb{R}^{m}$ such that
$y\notin\Phi_{1}^{+}(t;\mathbb{R}^{n}).$

In the first case,  $y$ is the image of a certain point from $\mathbb{R}^{n},$
namely the point
$x=\Phi_{1}(t)y,$ according to the equation $y=\Phi_{1}^{+}(t)x.$ In this case, we have $y=\Phi_{1}^{+}(t)\Phi_{1}(t)y, \, y=M_{1}(t)y,
\,  y\in M^{n}(t),$ which contradicts the assumption.

In the second case, we have $y=M_{1}(t)y$ and, according to (\ref{9}), $y=M_{1}(t)c(t)$
for a certain  $c(t)\in \mathbb{R}^{m}.$ Thus,
$y=M_{1}(t)c(t)=\Phi_{1}^{+}(t)\Phi_{1}(t)c(t)=\Phi_{1}^{+}(t)x,$ where
$x=\Phi_{1}(t)c(t),$ i.e.,  $y$ is the image of the point  $x$ under the mapping
$\Phi_{1}^{+}(t):\mathbb{R}^{n}\rightarrow \mathbb{R}^{m},$ whence $y\in
\Phi_{1}^{+}(t;\mathbb{R}^{n})$. This contradicts the assumption.

%%%%%%%%%%%%%%%%%%%%%%
\begin{center}\textbf{ 2. Main Theorem}\end{center}

\rm According to the results presented above, every continuously differentiable
nonsingular matrix
$\left(\begin{array}{cc}\Phi_{1}(t)\\\Phi_{2}(t)\end{array}\right)$ defines, in
the  $(t, y)$-space $\mathbb{R}\times \mathbb{R}^{m},  m=\dim
\left(\begin{array}{cc}\Phi_{1}(t)\\\Phi_{2}(t)\end{array}\right),$ the two
subspaces
\begin{equation}
\ker \,\Phi_{1}(t)=M_{1}^{n}(t),\,\,\,\,n=\dim\Phi_{1}(t), \nonumber
\end{equation}
\begin{equation} \ker \,\Phi_{2}(t)=M_{1}^{m-n}(t),\,\,\,\,p=\dim\Phi_{2}(t)=m-n,  \nonumber
\end{equation}
and two diffeomorphisms
\begin{equation}
\Phi_{1}^{+}(t):\mathbb{R}^{n}\rightarrow M_{1}^{n}(t),  \nonumber
\end{equation}
\begin{equation}
\Phi_{2}^{+}(t):\mathbb{R}^{p}\rightarrow M_{1}^{m-n}(t). \nonumber
\end{equation}
Using a matrix $Q(t)\in M_{m}(\mathbb{R})$ continuous for all  $t\in
\mathbb{R},$ we introduce, in the $(t, y)$-space $\mathbb{R}\times
\mathbb{R}^{m},$ a linear vector field  $(t,y'=Q(t)y)$ the integral curves of
which are defined by the solutions  $y=y(t)$ of the differential equation
 \setcounter{equation}{0}
\begin{equation}
\label{1} \frac{dy}{dt}=Q(t)y.
\end{equation}
If the union of the subspaces  $M_{1}^{n}(t)$ and $M_{1}^{m-n}(t)$ (or one
of these subspaces) is the union of integral curves of the vector field $(t, y'),$ then these
subspaces are called invariant manifolds of the differential equation (\ref{1})
or the vector field $(t, y').$ If the subspace
$M_{1}^{k}(t),\,\,\,k \in \{n, p\},$ is an invariant manifold of Eq.~(\ref{1}),
then  ``the motion of its points  $y$ in the  $(t, y)$-space is independent of
the motion of the points  $y$  outside the subspace  $M^{n}(t)$ for both $t>0$
and $t<0$.''

We pose the problem as follows: Find conditions under which the
subspace $M_{1}^{k}(t)$ is an invariant manifold of Eq.~(\ref{1}). An
equivalent statement of this problem is the following: Find conditions under which the
solutions $y=y(t) $ of Eq.~(\ref{1}) satisfy one of the additional conditions
\begin{equation}   \label{2}
y=M_{1}(t)y
\end{equation}
and
\begin{equation} \label{3}  M_{1}(t)y=0
\end{equation}
for any  $t\in \mathbb{R}.$

Finally, according to the terminology of Krylov--Bogolyubov nonlinear
mechanics [1, 2], the invariant manifold  $M^{k}(t)$ of Eq.~(\ref{1}) is an
integral manifold of Eq.~(\ref{1}) if, for any solution $y=y(t)$ of
Eq.~(\ref{1}), the fact that the inclusion
\begin{equation}
y(t)\in M^{k}(t),  \nonumber
\end{equation}
holds for a certain  $t=t_{0}$ implies that this inclusion is true for any $t\in
\mathbb{R}.$

Therefore, the posed problem is equivalent to the problem of finding
conditions under which the subspace $M^{k}(t)$ is an integral manifold of
Eq.~(\ref{1}).

\textbf{Theorem 1.} {\em Suppose that \/} $Q(t)\in M_{m}(\mathbb{R}),$ $\Phi(t)\in M_{m
n}(\mathbb{R}),$ $m>n,$ $\Phi^{+}(t)\in M_{m n}(\mathbb{R}),$ and  $ \rank\,
\Phi(t)=n.$ {\em Let \/} $Q(t)$ {\em be a continuous function and let
\/} $\Phi(t)$ {\em and \/} $\Phi^{+}(t)$ {\em be continuously
differentiable functions for all  \/} $t\in\mathbb{R}.$ {\em Also assume that \/}
$\Phi^{+}(t)$ {\em is a matrix pseudoinverse to the matrix \/} $\Phi(t)$ {\em
and \/}
\begin{gather*}
M^{n}(t)=\{y\in\mathbb{R}^{m}:\, y=M(t)y \},
\\
M^{m-n}(t)=\{y\in\mathbb{R}^{m}:\, M(t)y=0 \}
\\
M(t)=\Phi^{+}(t)\Phi(t),
\\
L(M,Q)=\frac{dM}{dt} + MQ-QM.
\end{gather*}

{\em Then the following assertions are true:\/}

1. {\em The subspaces \/} $M^{n}(t) $ {\em and \/} $M^{m-n}(t),$ {\em taken
together, are invariant manifolds of the differential equation \/}
$$
\frac{dy}{dt}=Q(t)y  \eqno(I)
$$
{\em if and only if \/}
$$
L(M(t), Q(t))=0.
$$

2. {\em The subspace  $M^{n}(t) $ is an invariant manifold of the differential
equation (I) if and only if  \/}
$$
L(M(t), Q(t)) M(t)=0
$$
{\em for any \/} $t\in \mathbb{R}. $ {\em Moreover, if \/} $M^{n}(t) $ {\em
is an invariant manifold of Eq.~(I), then, on \/} $M^{n}(t) $ {\em defined by
the diffeomorphism \/} $ \Phi^{+}(t),$
$$
y=\Phi^{+}(t)x, \qquad  x\in \mathbb{R}^{n},
$$
{\em Eq.~(I) is equivalent to the equation \/}
$$
\frac{dx}{dt}=P(t)x  \eqno(II)
$$
{\em with the coefficient matrix \/}
$$
P(t)=\left (\frac{d\Phi(t)}{dt}+\Phi(t)Q(t) \right)\Phi^{+}(t),
$$
{\em i.e., the fundamental matrices of solutions of Eqs.~(I) and (II)} $ Y(t)$ {\em and \/} $X(t),$ $Y(0)=E,$ $X(0)=E \in M_{n}(\mathbb{R}),$ {\em satisfy the relations
\/}
$$
Y(t)\Phi^{+}(0)=\Phi^{+}(t)X(t),
$$
$$
X(t)=\Phi(t)Y(t)\Phi^{+}(0)
$$
{\em for any  \/} $ t\in \mathbb{R}.$

3. {\em If \/} $ M^{m-n}(t)$ {\em is an invariant manifold of Eq.~(I), then \/}
$$
\ker \, L(M(t),Q(t))\supset M^{m-n}(t)
$$
{\em for any \/} $ t\in \mathbb{R}.$
\vspace*{2.0mm}  \\

We now pass to the proof of the theorem. Let
 \setcounter{equation}{3}
\begin{equation}   \label{4}
L(M(t),Q(t))=0 \qquad \forall \, t\in \mathbb{R}.
\end{equation}
Consider the function
\begin{equation}   \label{5}
r=(E-M(t))y(t),
\end{equation}
where $ y=y(t)$ is a solution of Eq.~(1) corresponding to the initial
conditions
\begin{equation}  \label{6}
y(t_{0})=M(t_{0})c
\end{equation}
and $c$ is an arbitrary point of the space $\mathbb{R}^m. $ According to
definition (\ref{5}), the function  $r$ is equal to zero for $t=t_0:$
\begin{equation}   \label{7}
r=r^0=(E-M(t_0))y(t_0)=(E-M(t_0))M(t_0)c=0.
\end{equation}
Differentiating function (\ref{5}), we obtain
\begin{equation}
\frac{dr}{dt}=-\frac{dM(t)}{dt}+(E-M(t))Q(t)y(t)= \nonumber
\end{equation}
\begin{equation}
=Q(t)y(t)-\left( \frac{dM(t)}{dt}+M(t)Q(t)-Q(t)M(t)\right)y(t)-Q(t)M(t)y(t)=
\nonumber
\end{equation}
\begin{equation}   \label{8}
=-L(M(t), Q(t)y(t)+Q(t)( y(t)-M(t)y(t))=Q(t)r.
\end{equation}
According to (\ref{7}), it follows from (\ref{8}) that $r(t)=0. $ Therefore,
\begin{equation}   \label{9}
y(t)=M(t)y(t) \qquad \forall \, t\in \mathbb{R}.
\end{equation}
On the one hand, we have
$$
\rank \, M(t)=\rank \,( \Phi^+ (t) \Phi(t)) \leq \min ( \rank
\, \Phi^+ (t),\ \rank \, \Phi(t))=n,
$$
while, on the other hand,
$$
n= \rank \,
(\Phi(t)\Phi^+ (t)\Phi(t)\Phi^+ (t))= \rank \, (\Phi(t)M(t)\Phi^+ (t))\leq
\rank \,M(t).
$$
Therefore,  $\rank \,M(t)=n  $ for any  $ t\in \mathbb{R}.$ Then
the subspace of  $ \mathbb{R}^m $ defined by points (\ref{5}) is
$n$-dimensional. Since points (\ref{6}) belong to the subspace  $ M^n (t_0), $
the subspace  $ M^n (t_0) $ coincides with the subspace defined by Eq.~(\ref{6}).
In this case, equality (\ref{9}) means that
\begin{equation}
y(t)\in M^n(t) \qquad \forall \, t \in \mathbb{R}  \nonumber
\end{equation}
for any solution of Eq.~(1) with initial value $y(t_0)=M(t_0)c $ for an
arbitrary  $c \in \mathbb{R}^m. $ Thus, the subspace  $M^n(t) $ is an invariant
manifold of Eq.~(1).

We now find the solution $ y=y(t) $ of Eq.~(1) with the initial conditions
\begin{equation}   \label{10}
y(t_0)=(E-M(t_0))c,
\end{equation}
where $ c $ is an arbitrary point of $ \mathbb{R}^m ,$ and consider the
function
\begin{equation}   \label{11}
r_1=M(t)y(t).
\end{equation}
Differentiating this function, we get
\begin{equation}
 \frac{dr_1}{dt}=\frac{dM(t)}{dt}y(t)+M(t)Q(t)y(t)=  \nonumber
\end{equation}
\begin{equation}=\left( \frac{dM(t)}{dt}+M(t)Q(t)-Q(t)M(t)\right)y(t)+Q(t)M(t)y(t)=
 \nonumber
\end{equation}
\begin{equation}   \label{12}
=L(M(t), Q(t))y(t)+Q(t)r=Q(t)r.
\end{equation}
By definition, the function $r_1$ is
equal to zero at the point $ t=t_0:$
\begin{equation}   \label{13}
r_1=r_1^0=M(t_0)y(t_0)=M(t_0)(E-M(t_0))c=0.
\end{equation}
Therefore, it follows from (\ref{12}) and (\ref{13}) that $r_1(t)=0 $ for any $
t \in \mathbb{R}. $ Thus,
\begin{equation}
 M(t)y(t)=0 \qquad \forall \, t \in \mathbb{R}, \nonumber
\end{equation}
which completes the proof of the inclusion
\begin{equation}   \label{13-1}
y(t)\in M^{m-n} (t)
\end{equation}
for any  $ t \in \mathbb{R}. $

Consider  $\rank\,(E-M(t_0))=m-\rank\, M(t_0)=m-n. $ Thus, the subspace
formed by points (\ref{10}) is $ (m-n)$-dimensional and coincides with the subspace
$ M^{m-n}(t_0)$. In this case, inclusion (\ref{13-1}) means that the subspace $
M^{m-n}(t)$ is an invariant manifold of Eq.~(1).

We have proved that the condition  $ L(M(t), Q(t))=0$ is sufficient for the
subspaces $ M^{n}(t)$ and $ M^{m-n}(t),$ taken together, to be invariant
manifolds of Eq.~(1).

Let the subspace  $ M^{n}(t)$ be an invariant manifold of Eq.~(1). Consider
the solutions of Eq.~(1)
\begin{equation}   \label{14}
y=Y(t) \Phi ^+ (0)c,
\end{equation}
where $ c $ is an arbitrary point of $ \mathbb{R}^m. $ The relation
\begin{equation} \label{15}
y(0)=Y(0) \Phi ^+ (0) c=\Phi ^+ (0) c \in M^{n} (0)
\end{equation}
yields the inclusion
\begin{equation}
Y(t) \Phi^{+}(0)c \in M^{n} \qquad \forall \, t \in \mathbb{R}.  \nonumber
\end{equation}
This means that
\begin{equation}   \label{16}
Y(t) \Phi^{+}(0)c = M(t)Y(t) \Phi^{+}(0)c
\end{equation}
for any  $ c \in \mathbb{R}^m $ and, hence, for unit vectors of the
space $\mathbb{R}^m .$ It follows from (\ref{16}) that
\begin{equation} \label{17}
Y(t) \Phi ^+ (0)=M(t)Y(t) \Phi^{+}(0)
\end{equation}
for any  $ t \in \mathbb{R}.$

Let $ X_t$ denote the matrix
\begin{equation} \label{18}
X_t=\Phi(t) Y(t) \Phi ^{+}(0).
\end{equation}
We rewrite  (\ref{17}) in the form of the relation
\begin{equation} \label{19}
Y(t) \Phi ^+ (0)=\Phi ^+ (t) X_t,
\end{equation}
which is true for any  $ t \in \mathbb{R}.$ Differentiating (\ref{19}) with
regard for (\ref{17}) and (\ref{18}), we obtain
\begin{equation} \label{20}
Q(t)Y(t) \Phi ^{+}(0)=\frac{d \Phi ^{+}(t)}{dt}X_t +\Phi ^{+}(t)
\frac{d X_t}{dt},
\end{equation}
\begin{equation}
Q(t)\Phi ^{+}(t)X_t=\frac{d \Phi ^{+}(t)}{dt}X_t +\Phi ^{+}(t)
\left( \frac{d \Phi(t)}{dt}X_t \Phi ^{+}(0)+ \Phi(t) Q(t)Y(t) \Phi
^{+}(0) \right)= \nonumber \end{equation}
\begin{equation} \label{21}
=\frac{d \Phi ^{+}(t)}{dt}X_t +\Phi ^{+}(t) \left( \frac{d
\Phi(t)}{dt}+\Phi(t) Q(t)\right) \Phi ^{+}(t)X_t.
\end{equation}
Subtracting (\ref{20}) from (\ref{21}), we get
\begin{equation}
\Phi ^{+}(t) \frac{d X_t}{dt}=\Phi ^{+}(t) \left( \frac{d \Phi(t)}{dt}+\Phi(t)
Q(t)\right) \Phi ^{+}(t)X_t.   \nonumber
\end{equation}
This proves that
\begin{equation}  \label{22}
 \Phi ^{+}(t) \left[ \frac{d X_t}{dt}- \left( \frac{d
\Phi(t)}{dt}+\Phi(t) Q(t)\right) \Phi ^{+}(t)X_t \right]=0
\end{equation}
for any  $ t \in \mathbb{R}.$ Since $ \ker \Phi ^{+}(t)=0,$ equality (\ref{22})
is possible only if
\begin{equation} \label{23}
\frac{d X_t}{dt}=P(t)X_t,
\end{equation}
where
\begin{equation} \label{24}
P(t)=\left( \frac{d \Phi(t)}{dt}+\Phi(t) Q(t)\right) \Phi ^{+}(t).
\end{equation}
By definition, we have
\begin{equation} \label{25}
X_0=\Phi(0)Y(0)\Phi ^{+}(0)=E, \qquad E \in M_n .
\end{equation}
Therefore,
\begin{equation} \label{26}
X_t=X(t)\qquad \forall \, t \in \mathbb{R}.
\end{equation}
Thus, if $ M^n (t)$ is an invariant manifold of Eq.~(1), then (\ref{19}) takes
the form
\begin{equation}   \label{27}
Y(t)\Phi^{+}(0)=\Phi^{+}(t)X(t),
\end{equation}
where $X(t),\,\,\,\,X(0)=E,$ is the fundamental matrix of solutions of
Eq.~(\ref{23}).

Multiplying (\ref{27}) from the left by  $\Phi(t),$ we get
\begin{equation}   \label{28}
X(t)=\Phi(t)Y(t)\Phi^{+}(0).
\end{equation}
Differentiating (\ref{17}), we obtain
\begin{equation}   \label{29}
QY(t)\Phi^{+}(0)=\left(\frac{dM(t)}{dt}+M(t)Q(t)-Q(t)M(t)\right)Y(t)\Phi^{+}(0)
+QM(t)Y(t)\Phi^{+}(0).
\end{equation}
In view of  (\ref{17}), relation (\ref{29}) yields
\begin{equation}   \label{30}
L(M(t),Q(t))Y(t)\Phi^{+}(0)=0.
\end{equation}
Using (\ref{30}) and (\ref{19}), we get
\begin{equation}
L(M(t),Q(t))\Phi^{+}(t)X(t)=0.   \nonumber
\end{equation}
Multiplying this equality from the right by $X^{-1}(t),$ we obtain the final
result
\begin{equation}   \label{31}
L(M(t),Q(t))\Phi^{+}(t)=0
\end{equation}
for any  $t\in \mathbb{R}.$

Thus, the fact that the subspace  $M^{n}(t)$ is integral implies that all
conditions of the theorem related to this case are satisfied; to this end, it suffices to rewrite
equality (\ref{27}) as an equation of the subspace  $M^{n}(t)$ in the
parametric form:
\begin{equation}  \label{32}
y=\Phi^{+}(t)x,\,\,\,\,x\in \mathbb{R}^{n},\,\,\,\,t\in
\mathbb{R}.
\end{equation}

Assume that the condition
\begin{equation}   \label{101}
L(M(t),Q(t))\Phi^{+}(t)=0 \qquad \forall \, t \in \mathbb{R},
\end{equation}
is satisfied. Multiplying  (\ref{101}) from the right by the matrix
$\Phi(t)X(t),$ where $X(t)$ is the fundamental matrix of solutions of
Eq.~(\ref{23}) with coefficient matrix (\ref{24}), $X(0)=E,$ we obtain
\begin{equation} \label{102}
L(M(t),Q(t))M(t)X(t)=0 \qquad \forall \, t \in \mathbb{R}.
\end{equation}
Now consider the function
\begin{equation}   \label{103}
r=y(t)-\Phi^{+}(t)X(t)c,
\end{equation}
where $y=y(t)$ is the solution of Eq.~(\ref{101}) such that
\begin{equation}   \label{104}
y(t_{0})=\Phi^{+}(t_{0})X(t_{0})c
\end{equation}
and $X(t),\,\,X(0)=E,$ is the fundamental matrix of solutions of the equation
\begin{equation}   \label{105}
\frac{dx}{dt}=P(t)x
\end{equation}
with coefficient matrix (\ref{24}). Differentiating function (\ref{103}),
we obtain
\begin{equation}
\frac{dr}{dt}=Q(t)y(t)-\left(\frac{d\Phi^{+}(t)}{dt}X(t)c+\Phi^{+}(t)P(t)X(t)c\right)=
Q(t)[\,(y(t)-\Phi^{+}(t)X(t)c)\,]+  \nonumber
\end{equation}
\begin{equation}
+Q(t)\Phi^{+}(t)X(t)c-\left(\frac{d\Phi^{+}(t)}{dt}X(t)c+\Phi^{+}(t)P(t)X(t)c\right)=
 \nonumber
\end{equation}
\begin{equation}   \label{106}
=Q(t)r-\left(\frac{d\Phi^{+}(t)}{dt}+\Phi^{+}(t)P(t)-Q(t)\Phi^{+}(t)\right)X(t)c\,\,.
\end{equation}
Let us prove that the second term in (\ref{106}) is $0\in M_{mn}.$ Indeed, taking
(\ref{24}) into account, we get
\begin{equation}
\frac{d\Phi^{+}(t)}{dt}+\Phi^{+}(t)\left(\frac{d\Phi(t)}{dt}+\Phi(t)Q(t)\right)\Phi^{+}(t)
-Q(t)\Phi^{+}(t)=  \nonumber
\end{equation}
\begin{equation}
=\frac{d\Phi^{+}(t)}{dt}
+\Phi^{+}(t)\frac{d\Phi(t)}{dt}\Phi^{+}(t)+M(t)Q(t)\Phi^{+}(t)-Q(t)\Phi^{+}(t)=
\nonumber
\end{equation}
\begin{equation}
=\left(\frac{d\Phi^{+}(t)}{dt}\Phi(t)+\Phi^{+}(t)\frac{d\Phi(t)}{dt}\right)\Phi^{+}(t),\nonumber
\end{equation}
\begin{equation}
M(t)Q(t)\Phi^{+}(t)-Q(t)M(t)\Phi^{+}(t)+Q(t)M(t)\Phi^{+}(t)-Q(t)\Phi^{+}(t)=\nonumber
\end{equation}
\begin{equation}   \label{107}
=\left(\frac{dM(t)}{dt}+M(t)Q(t)-Q(t)M(t)\right)\Phi^{+}(t)+
\end{equation}
\begin{equation}
+Q(t)\Phi^{+}(t)\Phi(t)\Phi^{+}(t)-Q(t)\Phi^{+}(t)=0.  \nonumber
\end{equation}
With regard for (\ref{107}), equality (\ref{106}) takes the form
\begin{equation}   \label{108}
\frac{dr}{dt}=Q(t)r.
\end{equation}
For $t=t_{0},$ according to (\ref{104}), function (\ref{103}) is equal
to zero:
\begin{equation}
r(t_{0})=y(t_{0})-\Phi^{+}(t_{0})X(t_{0})c=0.  \nonumber
\end{equation}
Therefore, it follows from (\ref{108}) that
\begin{equation}   \label{109}
r(t)=0 \qquad \forall \, t \in \mathbb{R},
\end{equation}
and, hence,
\begin{equation}   \label{110}
y(t)=\Phi^{+}(t)X(t)c \qquad \forall \, t \in \mathbb{R},
\end{equation}
where $c$ is an arbitrary point of the space $\mathbb{R}^{n}.$

Since (\ref{110}) is the parametric representation of the equation of the subspace
$M^{n}(t),$ it follows from (\ref{110}) that the condition that
\begin{equation}   \label{111}
y(t)\in M^{n}(t)
\end{equation}
for $t=t_{0}$ implies that inclusion (\ref{111}) holds for any $t\in
\mathbb{R}$. This proves that condition (\ref{101}) is not only necessary but
also sufficient for the subspace  $M^{n}(t)$ to be an invariant
manifold of Eq.~(\ref{1}).

This completes the proof of the second assertion of Theorem~1.

Let the subspace  $M^{m-n}(t)$ be an invariant manifold of Eq.~(\ref{1}). Consider the solutions  $y=y(t)$ of Eq.~(\ref{1}) defined by the relation
\begin{equation}
y(t)=Y(t)(E-M(0))c,  \nonumber
\end{equation}
where $c$ is an arbitrary point of  $\mathbb{R}^{m}.$

It follows from the relation
\begin{equation}
(0)y(0)=M(0)Y(0)(E-M(0))c=M(0)(E-M(0))c=0 \nonumber
\end{equation}
that  $y(0)\in M^{m-n}(0)$ and, hence,  $ y(t)\in M^{m-n}(t)
\quad \forall \, t \in \mathbb{R}.$ This proves that
\begin{equation}   \label{32'}
M(t)Y(t)(E-M(0))c=0 \qquad \forall \, t \in \mathbb{R}.
\end{equation}
Differentiating this equality, we get
\begin{equation}
 \left(\frac{dM(t)}{dt}+M(t)Q(t)-Q(t)M(t)\right)Y(t)(E-M(0))c\,+\nonumber
\end{equation}
\begin{equation}
\label{33}+Q(t)M(t)Y(t)(E-M(0))c=L(M(t),Q(t))Y(t)(E-M(0))c=0 .
\end{equation}
The points  $y=(E-M(0))c$ define the subspace $M^{m-n}(0).$ Therefore, the
equation
\begin{equation}   \label{34}
y=Y(t)(E-M(0))c, \,\,\,\,c\in \mathbb{R}^{m},\,\,\,\,t\in \mathbb{R},
\end{equation}
defines the subspace  $M^{m-n}(t)$ in the parametric form. Therefore, equality
(\ref{33}) means that
\begin{equation}   \label{35}
\ker \,L(M(t),Q(t))\supset M^{m-n}(t) \qquad \forall \, t
\in \mathbb{R}.
\end{equation}
Thus, inclusion (\ref{35}) is a necessary condition for the subspace
$M^{m-n}(t)$ to be an invariant manifold of Eq.~(\ref{1}).

Assume that the subspaces $M^{n}(t)$ and $M^{m-n}(t),$ taken together, are invariant
manifolds of Eq.~(\ref{1}). Then equality (\ref{30}) yields
\begin{equation}   \label{36}
\ker\,L(M(t),Q(t))\supset Y(t)M(0)c \qquad \forall \, t \in \mathbb{R},
\end{equation}
where $c$ is an arbitrary point of  $\mathbb{R}^{m}.$ Moreover, since the
equation
\begin{equation}
y=Y(t)M(0)c, \,\,\,c\in \mathbb{\mathbb{R}}^{m},\,\,\,\,t\in \mathbb{R}, \nonumber
\end{equation}
defines the subspace  $M^{n}(t)$ in the parametric form, it follows from (\ref{36}) that
\begin{equation}   \label{37}
\ker \,L(M(t),Q(t))\supset M^{n}(t) \qquad \forall \, t \in \mathbb{R}.
\end{equation}
Thus, if the subspaces  $M^{n}(t)$ and $M^{m-n}(t),$ taken together, are
invariant manifolds of Eq.~(\ref{1}), then
\begin{equation}   \label{38}
\ker \,L(M(t),Q(t))\supset (M^{n}(t)\cup M^{m-n}(t)) .
\end{equation}
Since $\rank\,\,M^{n}(t)=n,\,\,\,\,\rank \,M^{m-n}(t)=m-n,$ and $M^{n}(t)\cap
M^{m-n}(t)=\{0\},$ we conclude that the union on the right-hand side of expression (\ref{38})
contains a basis of the space  $\mathbb{R}^{m}.$ Therefore, relation (\ref{38}) yields
\begin{equation}   \label{39}
\ker \,L(M(t),Q(t))\supset \mathbb{R}^{m} \qquad \forall
\, t \in \mathbb{R}.
\end{equation}
Since $ \ker \,L(M(t),Q(t))\in M_m (\mathbb{R}) \quad \forall \, t \in
\mathbb{R},$ inclusion (\ref{39}) is possible only if
\begin{equation}   \label{40}
\ker \,L(M(t),Q(t))=0.
\end{equation}
Thus, condition (\ref{40}) is not only sufficient but also necessary for the
subspaces $M^{n}(t)$ and $M^{m-n}(t),$ taken together, to be invariant manifolds
of Eq.~(\ref{1}).

\vspace*{6.0mm}

%%%%%%

\noindent 1. N. N. Bogolyubov, {\it On Some Statistical Methods in
Mathematical Physics} [in Russian], Academy of Sciences of Ukr.~SSR, Kiev (1945); N. N.
Bogolyubov, {\it Collection of Scientific Works} [in Russian], Vol.~4, Nauka,
Moscow (2006).

\noindent 2.  N. N. Bogolyubov and Yu. A. Mitropol'skii, {\it Asymptotic
Methods in the Theory of Nonlinear Oscillations} [in Russian], Fizmatgiz,
Moscow (1963); N. N. Bogolyubov, {\it Collection of Scientific Works} [in
Russian], Vol.~3, Nauka, Moscow (2005).

\end{document}

%\fontsize{14.4}{18pt}\selectfont  % - для распечатки в 14-м кегле и
% через 2 интервала (как на печатной машинке через 1,5 интервала).
% чтоб распечатать через 1,5 интервала, вместо {24pt} поставить {18pt}.

УДК 517.9 \vspace{6.0mm}

 \textbf{А. М. Самойленко}\quad
(Ін-т математики НАН України, Київ)\vspace{6.0mm}

\textbf{ Про інваріанті многовиди лінійних диференціальних рівнянь.}
\vspace{6.0mm}

В роботі розвинено ідеї методу інтегральних многовидів М. М.
Боголюбова стосовно лінійних диференціальних рівнянь. Отриманий
результат може знайти прикладну інтерпритацію та застосування.
\vspace{6.0mm}
\begin{center}\textbf{ 1. Дві леми.}\end{center}

\rm Нехай матриці $\Phi_{1}(t)$ та $\Phi_{2}(t)$ неперервно
диференційовані для всіх  $t\in R$, $\Phi_{1}(t)\in M_{nm}(R)$,
$\Phi_{2}(t)\in M_{pm}(R)$, де $m>n>p=m-n$ та задовольняють умову
$$det\left(\begin{array}{cc}\Phi_{1}(t)\\\Phi_{2}(t)\end{array}\right)\ne 0\,\,\,\,.$$

Позначимо через $\Phi_{1}^{+}(t)$ та $\Phi_{2}^{+}(t)$
псевдообернені до $\Phi_{1}(t)$ та $\Phi_{2}(t)$ матриці, визначені
умовою
\begin{equation}
\label{1}
\left(\begin{array}{cc}\Phi_{1}(t)\\\Phi_{2}(t)\end{array}\right)(\Phi_{1}^{+}(t),\Phi_{2}^{+}(t))=E,
\,\,\,\,E\in M_{m}(R).
\end{equation}
Нехай
\begin{equation}
\label{2} M_{1}(t)=\Phi_{1}^{+}(t)\Phi_{1}(t),
\end{equation}
\begin{equation}
\label{3} M_{2}(t)=\Phi_{2}^{+}(t)\Phi_{2}(t).
\end{equation}
Із (~\ref{1}) та визначень матриць $M_{1}(t)$ та $M_{2}(t)$
випливає, що ці матриці задовольняють умови
\begin{equation}
M_{\nu}^{2}(t)=M_{1}(t),\,\,\,\,\nu=1,2,\,\,\,\,\rank\,M_{1}(t)=n,\,\,\,\,\rank\,M_{2}(t)=p,\,\,\,\,M_{1}(t)M_{2}(t)=M_{2}(t)M_{1}(t)=0,\nonumber
\end{equation}
а також
\begin{equation}
\label{4} M_{1}(t)+M_{2}(t)=E.
\end{equation}
Доведемо рівність (~\ref{4}). Справді, оскільки матриця
$(\Phi_{1}^{+}(t),\Phi_{2}^{+}(t))$ обернена до матриці\\\\
$\left(\begin{array}{cc}\Phi_{1}(t)\\\Phi_{2}(t)\end{array}\right)$,
то, домноживши останню зліва на матрицю
$(\Phi_{1}^{+}(t),\Phi_{2}^{+}(t))$, отримаємо рівність
$$\Phi_{1}^{+}(t)\Phi_{1}(t)+\Phi_{2}^{+}(t)\Phi_{2}(t)=E,$$
яка тотожна до рівності (\ref{4}).

Визначимо в просторі $R^{m}$ два підпростори, використавши матрицю
$M_{1}(t)$:
\begin{equation}
\label{5} M^{n}(t)=\{y\in R^{m}: y=M_{1}(t)\},
\end{equation}
\begin{equation}
\label{6} M^{m-n}(t)=\{y\in R^{m}: M_{1}(t)y=0\},
\end{equation}
та два підпростори, використавши матрицю $M_{2}(t)$:
\begin{equation}
\label{7} M_{1}^{p}(t)=\{y\in R^{m}: y=M_{2}(t)\},
\end{equation}
\begin{equation}
\label{8} M_{1}^{m-p}(t)=\{y\in R^{m}: M_{2}(t)y=0\}.
\end{equation}
%%%%%%%%%%% {\em и \/}

\textbf{Лема 1.}\,\,\,  {\em Підпростори \/} $M^{k}(t)$ {\em та \/}
$M_{1}^{k}(t),\,\,\,\,k\in\{n,p\}$, {\em задовольняють умови \/}
$$M^{n}(t)=M_{1}^{m-p}(t)=\ker \,\Phi_{2}(t),\,\,\,\,M^{m-n}(t)=M_{1}^{p}(t)=\ker \,\Phi_{1}(t).$$

Справді, за властивостями матриці $M_{1}(t)$ загальним розв'язком
рівняння, що визначає підпростір (~\ref{5}), є функція
\begin{equation}
\label{9} y=M_{1}(t)c(t),
\end{equation}
де $c(t)$ --- довільна функція із значеннями в $R^{m}$. За
властивостями матриці $M_{2}(t)$ загальним розв'язком рівняння, що
визначає підпростір (~\ref{8}), є функція
\begin{equation}
\label{10} y=(E-M_{2}(t))c_{1}(t)=M_{1}(t)c_{1}(t).
\end{equation}
Загальні розв'язки розглядуваних рівнянь співпадають при
$c(t)=c_{1}(t)$, що доводить рівність $M^{n}(t)=M_{1}^{m-p}(t)$.

Із визначення $M_{1}^{m-p}$ випливає, що
$$M_{1}^{m-p}(t)=\ker \,M_{2}(t)=\ker \,(\Phi_{2}^{+}(t)\Phi_{2}(t)).$$
А оскільки $\ker \,\Phi_{2}^{+}(t)=0$, то для кожного $x$, для якого
$\Phi_{2}^{+}(t)\Phi_{2}(t)x=0$, буде виконуватись співвідношення
$\Phi_{2}(t)x=0$, тому $\ker \,(\Phi_{2}^{+}(t)\Phi_{2}(t))=\ker
\,\Phi_{2}(t)$. Це доводить першу із рівностей леми 1. Друга її
рівність доводиться аналогічно.
%%%%%%%%%%% {\em и \/}

\textbf{Лема 2.}\,\,\, {\em Відображення \/}
$\Phi_{1}^{+}(t):\,\,\,\,y=\Phi_{1}^{+}(t)x$ {\em є
дифеоморфізмом\/} $R^{n}$ {\em в\/} $M^{n}(t)$, {\em а відображення
\/} $\Phi_{2}^{+}(t):\,\,\,\,y=\Phi_{2}^{+}(t)x$ {\em є
дифеоморфізмом\/} $R^{n}$ {\em в\/} $M^{m-n}(t)$.

Доведемо перше з тверджень леми 2, бо друге доводиться
аналогічно.\\

Матриця $\Phi_{1}^{+}(t)$, як блок матриці, оберненої до матриці
$\left(\begin{array}{cc}\Phi_{1}(t)\\\Phi_{2}(t)\end{array}\right)$,
є неперервно диференційованою та має неперервно диференційовану
псевдообернену матрицю, якою є матриця $\Phi_{1}(t)$. Тому для
доведення першого з тверджень леми 2 залишається довести, що
відображення $\Phi_{1}^{+}(t):R^{n}\rightarrow M^{n}(t)$ є
гомеоморфізмом.

Оскільки $\ker \Phi_{1}^{+}(t)=0$, то $\Phi_{1}^{+}(t)$ є
гомеоморфізмом $R^{n}$ в образ $\Phi_{1}^{+}(t;R^{n})$ при
відображенні $\Phi_{1}^{+}(t):R^{n}\rightarrow R^{m}$. Залишилося
довести, що $\Phi_{1}^{+}(t;R^{n})=M^{n}(t)$.

Нехай це не так. Тоді або існує така точка
$y\in\Phi_{1}^{+}(t;R^{n})$, що $y\notin M^{n}(t)$, або існує така
точка $y\in M^{n}(t)$, що $y\notin\Phi_{1}^{+}(t;R^{n})$.

В першому випадку існує точка $x\in R^{n}$ така, що
$y=\Phi_{1}^{+}(t)x$. Тоді $x=\Phi_{1}(t)$, а тому
$y=\Phi_{1}^{+}(t)\Phi_{1}(t)y=M_{1}(t)y$. Це означає, що $y\in
M^{n}(t)$, що суперечить припущенню.

В другому випадку $y=M_{1}(t)y$, тому існує таке $c(t)\in R^{m}$, що
$y=M_{1}(t)c(t)$. Отже за формулою (~\ref{9})
$y=M_{1}(t)c(t)=\Phi_{1}^{+}(t)\Phi_{1}(t)c(t)=\Phi_{1}^{+}(t)x$, де
$x=\Phi_{1}(t)c(t)$, тобто $y$ є образом точки $x$ при відображенні
$\Phi_{1}^{+}(t):R^{n}\rightarrow R^{m}$ , звідки $y\in
\Phi_{1}^{+}(t;R^{n})$. Це суперечить припущенню.

%%%%%%%%%%%%%%%%%%%%%%
\begin{center}\textbf{ 2. Основна теорема.}\end{center}

\rm Згідно з викладеним раніше кожна неперервно диференційована
матриця \\\\
$\left(\begin{array}{cc}\Phi_{1}(t)\\\Phi_{2}(t)\end{array}\right)$
у випадку її невиродженості визначає в $(t;y)$ просторі $R\times
R^{m},\,\,\,\,m=\dim \left(\begin{array}{cc}\Phi_{1}(t)\\\Phi_{2}(t)\end{array}\right)$\\\\
два підпростори
$$\ker \,\Phi_{1}(t)=M_{1}^{n}(t),\,\,\,\,n=\dim\Phi_{1}(t)$$
$$\ker \,\Phi_{2}(t)=M_{1}^{m-n}(t),\,\,\,\,p=\dim\Phi_{2}(t)=m-n$$
та два дифеоморфізми
$$\Phi_{1}^{+}(t):R^{n}\rightarrow M_{1}^{n}(t)$$
$$\Phi_{2}^{+}(t):R^{p}\rightarrow M_{1}^{m-n}(t).$$
Нехай $Q(t)\in M_{m}(R)$ --- матриця, неперервна для всіх $t\in R$ в
$(t;y)$ просторі $R\times R^{m}$. Розглянемо лінійне векторне поле
$(t,y'=Q(t)y)$, інтегральні криві якого визначаються розв'язком
$y=y(t)$ диференціального рівняння
\setcounter{equation}{0}
\begin{equation}
\label{1} \frac{dy}{dt}=Q(t)y).
\end{equation}
Якщо один з підпросторів $M_{1}^{n}(t)$ та $M_{1}^{m-n}(t)$ є
об'єднанням інтегральних кривих векторного поля $(t,y')$, то його
називають інваріантним підпростором диференціального \\рівняння
(\ref{1}) або векторного поля $(t,y')$. Коли підпростір
$M_{1}^{k}(t),\,\,\,\,k\in {n,p}$ є інваріантним відносно рівняння
(\ref{1}), то рух його точок $y$ в $(t;y)$ просторі відбувається
незалежно від того, як рухаються точки поза підпростором
$M_{1}^{k}(t)$ як для $t>0$, так і для $t<0$.

Поставимо задачу: знайти умови, при яких підпростір $M_{1}^{k}(t)$ є
інваріантним многовидом рівняння (\ref{1}). Еквівалентна постановка
цієї задачі: знайти умови, при яких розв'язки рівняння (\ref{1})
задовольняють одну з додаткових умов
\begin{equation}
\label{2} y=M_{1}(t)y
\end{equation}
або
\begin{equation} \label{3} M_{1}(t)y=0
\end{equation}
для будь-якого $t\in R$.

Нарешті, згідно з термінологією нелінійної механіки Крилова --
Боголюбова [\,1, 2\,] інваріантний многовид $M^{k}(t)$ рівняння
(\ref{1}) є інтегральним многовидом (\ref{1}), якщо для довільного
розв'язку $y=y(t)$ рівняння (\ref{1}) із співвідношення
$$y(t)\in M^{k}(t),$$
справедливого для якогось $t=t_{0}$, випливає його справедливість
для будь-якого $t\in R$.

Отже, поставлена нами задача також еквівалентна задачі про
знаходження умов, при яких підпростір $M^{k}(t)$ є інтегральним
многовидом рівняння (\ref{1}).

%%%%%%%%%%% {\em и\/}
\textbf{Теорема 1.} {\em Нехай \/} $Q(t)\in M_{m}(\mathbb{R}),\,
\Phi(t)\in M_{m n}(\mathbb{R}), \, m>n, \, \Phi^{+}(t)\in M_{m
n}(\mathbb{R}), \,  \rank\, \Phi(t)=n.$ {\em Крім того, нехай \/}
$Q(t)$ {\em є неперервною,\/} $\Phi(t)$ {\em та \/} $\Phi^{+}(t)$
{\em є неперервно диференційованими функціями для всіх \/}
$t\in\mathbb{R}$. {\em Більше того, нехай \/} $\Phi^{+}(t)$ {\em є
матрицею, псевдооберненою для матриці \/} $\Phi(t)$. {\em Нарешті,
нехай \/}
$$M^{n}(t)={y\in\mathbb{R}^{m}:\, y=M(t)y },$$
$$M^{m-n}(t)={y\in\mathbb{R}^{m}:\, M(t)y=0 }$$
{\em та \/}
$$M(t)=\Phi^{+}(t)\Phi(t),$$
$$L(M,Q)=\frac{dM}{dt} + MQ-QM. $$

{\em Тоді:\/}

1. {\em Підпростори \/} $M^{n}(t) $ {\em та \/} $M^{m-n}(t) $, {\em
взяті разом, є інваріантними многовидами диференціального рівняння
\/}
$$\frac{dy}{dt}=Q(t)y  \eqno(I)$$
{\em тоді і лише тоді, коли виконується умова \/}
$$ L(M(t), Q(t))=0. $$

2. {\em Підпростір  є інваріантним многовидом диференціального
рівняння (I) тоді і лише тоді, коли \/}
$$ L(M(t), Q(t)) M(t)=0 $$
{\em для будь-якого\/} $t\in \mathbb{R}. $ {\em Крім того, коли\/}
$M^{n}(t) $ {\em є інваріантним многовидом рівняння (I),  то на \/}
$M^{n}(t), $ {\em визначеним дифеоморфізмом \/} $ \Phi^{+}(t): $
$$ y=\Phi^{+}(t)x, \qquad  x\in \mathbb{R}^{n},$$
{\em рівняння (I) еквівалентне рівнянню \/}
$$\frac{dx}{dt}=P(t)x  \eqno(II)$$
{\em з матрицею коефіцієнтів \/}
$$P(t)=\left (\frac{d\Phi(t)}{dt}+\Phi(t)Q(t) \right)\Phi^{+}(t), $$
{\em тобто, фундаментальні матриці розв'язків рівнянь (I) та (II)\/}
$ Y(t)$ {\em і \/} $X(t), \; Y(0)=E$ {\em і \/} $X(0)=E \in
M_{n}(\mathbb{R})$ {\em задовольняють співвідношенням \/}
$$Y(t)\Phi^{+}(0)=\Phi^{+}(t)X(t), $$
$$X(t)=\Phi(t)Y(t)\Phi^{+}(0) $$
{\em для будь-якого \/} $ t\in \mathbb{R}.$

3. {\em Якщо \/} $ M^{m-n}(t)$ {\em є інваріантним многовидом
рівняння (I), то \/}
$$\ker \, L(M(t),Q(t))\supset M^{m-n}(t) $$
{\em для будь-якого \/} $ t\in \mathbb{R}.$
\vspace*{4.0mm}  \\
Перейдемо до доведення теореми. Нехай \setcounter{equation}{3}
\begin{equation}
\label{4} L(M(t),Q(t))=0 \qquad \forall \, t\in \mathbb{R}.
\end{equation}
Розглянемо функцію
\begin{equation}
\label{5} r=(E-M(t))y(t),
\end{equation}
де $ y=y(t)$ --- розв'язок рівняння (1), визначений початковими
умовами
\begin{equation}
\label{6} y(t_{0})=M(t_{0})c
\end{equation}
де $c$ --- довільна точка простору $\mathbb{R}^m. $ Згідно із
визначенням (\ref{5}), функція $r$ набуває при $t=t_0 $ нульового
значення
\begin{equation}
\label{7} r=r^0=(E-M(t_0))y(t_0)=(E-M(t_0))M(t_0)c=0.
\end{equation}
Диференціюючи функцію (\ref{5}), отримаємо, що
\begin{equation}
\frac{dz}{dt}=-\frac{dM(t)}{dt}+(E-M(t))Q(t)y(t)= \nonumber
\end{equation}
\begin{equation}=Q(t)y(t)-\left( \frac{dM(t)}{dt}+M(t)Q(t)-Q(t)M(t)\right)y(t)-Q(t)M(t)y(t)=
\nonumber \end{equation}
\begin{equation}
\label{8} =-L(M(t), Q(t)y(t)+Q(t)( y(t)-M(t)y(t))=Q(t)r.
\end{equation}
Згідно з (\ref{7}), із (\ref{8}) випливає, що $\varepsilon (t)=0, $
отже,
\begin{equation}
\label{9} y(t)=M(t)y(t) \qquad \forall \, t\in \mathbb{R}.
\end{equation}
Оскільки, з одного боку, $\rank \, M(t)=\rank \,( \Phi^+
(t)+\Phi(t)) \leq  \min ( \rank \, \Phi^+ (t),\, \rank \,
\Phi(t))=n,$ а з другого --- $n= \rank \, (\Phi(t)\Phi^+
(t)\Phi(t)\Phi^+ (t))= \rank \, (\Phi(t)M(t)\Phi^+ (t))\leq \rank
\,M(t),$ то $\rank \,M(t)=n  $ для будь-якого $ t\in \mathbb{R}.$
Тоді підпростір в $ \mathbb{R}^m ,$ що визначається точками
(\ref{5}), є $n$-мірним, і оскільки точки (\ref{6}) належать
підпростору $ M^n (t_0), $ то півпростір $ M^n (t_0) $ збігається з
підпростором, визначеним рівнянням (\ref{6}). Рівність (\ref{9})
означає тоді, що
$$ y(t)\in M^n(t) \qquad \forall \, t \in \mathbb{R} $$
для будь-якого розв'язку рівняння (1) з початковим значенням
$y(t_0)=M(t_0)c  $ для довільного $c \in \mathbb{R}^m. $ Це означає,
що підпростір $M^n(t) $ є інваріантним многовидом рівняння (1).

Визначимо тепер розв'язок $ y=y(t) $ рівняння (1) початковими
умовами
\begin{equation}
\label{10} y(t_0)=(E-M(t_0))c,
\end{equation}
де $ c $ --- довільна точка $ \mathbb{R}^m ,$ та розглянемо функцію
\begin{equation}
\label{11} r_1=M(t)y(t).
\end{equation}
Диференціюючи цю функцію, отримаємо, що
\begin{equation}
 \frac{dr_1}{dt}=\frac{dM(t)}{dt}y(t)+M(t)Q(t)= \nonumber\end{equation}
\begin{equation}=\left( \frac{dM(t)}{dt}+M(t)Q(t)-Q(t)M(t)\right)y(t)+Q(t)M(t)y(t)=
\nonumber \end{equation}
\begin{equation}
\label{12}=L(M(t), Q(t))y(t)+Q(t)r=Q(t)r.
\end{equation}
Оскільки згідно з визначенням функція $r_1 $ в точці $ t=t_0 $
набуває нульвого значення
\begin{equation}
\label{13} r_1=r_1^0=M(t_0)y(t_0)=M(t_0)(E-M(t_0))c=0,
\end{equation}
то із (\ref{12}) та (\ref{13}) випливає, що $r_1(t)=0 $ для
будь-якого $ t \in \mathbb{R}, $ отже,
\begin{equation}
 M(t)y(t)=0 \qquad \forall \, t \in \mathbb{R}, \nonumber
\end{equation} що доводить справедливість включення
\begin{equation}
\label{13-1} y(t)\in M^{m-n} (t)
\end{equation}
для будь-якого $ t \in \mathbb{R}. $

Розглянемо $\rank\,(E-M(t_0))=m-\rank\, M(t_0)=m-n. $ Отже,
підпростір, створений точками (\ref{10}), є $ (m-n)$-вимірним та
збігається з підпростором $ M^{m-n}(t_0)$. Включення (\ref{13-1})
тоді означає, що підпростір $ M^{m-n}(t)$ є інваріантним многовидом
рівняння (1).

Нами доведено, що умова $ L(M(t), Q(t))=0$ є достатньою для того,
щоб підпростори $ M^{n}(t)$ та $ M^{m-n}(t),$ взяті разом, були
інваріантними многовидами рівняння (1).

Нехай підпростір $ M^{n}(t)$ є інваріантним многовидом рівняння (1).
Розглянемо розв'язки рівняння (1)
\begin{equation} \label{14}
y=Y(t) \Phi ^+ (0)c,
\end{equation}
де $ c $ --- довільна точка $ \mathbb{R}^m. $ Із того, що
\begin{equation} \label{15}
y(0)=Y(0) \Phi ^+ (0) c=\Phi ^+ (0) c \in M^{n} (0),
\end{equation}
випливає включення
\begin{equation}  Y(t) \Phi^{+}(0)c \in M^{n} \qquad \forall \, t \in
\mathbb{R}.
\nonumber \end{equation}
Це означає, що
\begin{equation}
\label{16} Y(t) \Phi^{+}(0)c = M(t)Y(t) \Phi^{+}(0)c
\end{equation} для довільного
$ c \in \mathbb{R}^m ,$ отже, і для одиничних ортів простору
$\mathbb{R}^m .$ Тоді із (\ref{16}) випливає співвідношення
\begin{equation} \label{17}
Y(t) \Phi ^+ (0)=M(t)Y(t) \Phi^{+}(0)
\end{equation}
для будь-якого $ t \in \mathbb{R}.$

Позначимо через $ X_t$ матрицю
\begin{equation} \label{18}
X_t=\Phi(t) Y(t) \Phi ^{+}(0)
\end{equation}
та запишемо (\ref{17}) у вигляді співвідношення
\begin{equation} \label{19}
Y(t) \Phi ^+ (0)=\Phi ^+ (t) X_t,
\end{equation}
яке справджується для будь-якого $ t \in \mathbb{R}.$ Диференціюючи
(\ref{19}),
 з урахуванням співвідношень (\ref{17}), (\ref{18})
отримуємо, що
\begin{equation} \label{20}
Q(t)Y(t) \Phi ^{+}(0)=\frac{d \Phi ^{+}(t)}{dt}X_t +\Phi ^{+}(t)
\frac{d X_t}{dt},
\end{equation}
\begin{equation}
Q(t)\Phi ^{+}(t)X_t=\frac{d \Phi ^{+}(t)}{dt}X_t +\Phi ^{+}(t)
\left( \frac{d \Phi(t)}{dt}X_t \Phi ^{+}(0)+ \Phi(t) Q(t)Y(t) \Phi
^{+}(0) \right)= \nonumber \end{equation}
\begin{equation} \label{21}
=\frac{d \Phi ^{+}(t)}{dt}X_t +\Phi ^{+}(t) \left( \frac{d
\Phi(t)}{dt}+\Phi(t) Q(t)\right) \Phi ^{+}(t)X_t.
\end{equation}
Віднявши (\ref{20}) від (\ref{21}), отримаємо рівність
\begin{equation}
\Phi ^{+}(t) \frac{d X_t}{dt}=\Phi ^{+}(t) \left( \frac{d
\Phi(t)}{dt}+\Phi(t) Q(t)\right) \Phi ^{+}(t)X_t.
 \nonumber \end{equation} Це доводить,
що
\begin{equation}
 \label{22}
 \Phi ^{+}(t) \left[ \frac{d X_t}{dt}- \left( \frac{d
\Phi(t)}{dt}+\Phi(t) Q(t)\right) \Phi ^{+}(t)X_t \right]=0
\end{equation}
для будь-якого $ t \in \mathbb{R}.$ Оскільки $ \ker \Phi ^{+}(t)=0,$
то рівність (\ref{22}) можлива лише тоді, коли
\begin{equation} \label{23}
\frac{d X_t}{dt}=P(t)X_t,
\end{equation}
де
\begin{equation} \label{24}
P(t)=\left( \frac{d \Phi(t)}{dt}+\Phi(t) Q(t)\right) \Phi ^{+}(t).
\end{equation}
Оскільки за визначенням
\begin{equation} \label{25}
X_0=\Phi(0)Y(0)\Phi ^{+}(0)=E, \qquad E \in M_n ,
\end{equation}
то \begin{equation} \label{26} X_t=X(t)\qquad \forall \, t \in
\mathbb{R}.
\end{equation}
Отже, коли $ M^n (t)$ є інваріантним многовидом рівняння (1), то
(\ref{19}) набуває вигляду
\begin{equation}
\label{27}Y(t)\Phi^{+}(0)=\Phi^{+}(t)X(t),
\end{equation}
де $X(t),\,\,\,\,X(0)=E$ --- фундаментальна матриця розв'язків
рівняння (~\ref{23}).

Домножимо (~\ref{27}) зліва на $\Phi(t)$ та отримаємо
\begin{equation}
\label{28}X(t)=\Phi(t)Y(t)\Phi^{+}(0).
\end{equation}
Продиференціювавши (~\ref{17}), одержимо
\begin{equation}
\label{29}QY(t)\Phi^{+}(0)=\left(\frac{dM(t)}{dt}+M(t)Q(t)-Q(t)M(t)\right)Y(t)\Phi^{+}(0)
+QM(t)Y(t)\Phi^{+}(0).
\end{equation}
Із (~\ref{29}) з урахуванням (~\ref{17}) випливає, що
\begin{equation}
\label{30} L(M(t),Q(t))Y(t)\Phi^{+}(0)=0.
\end{equation}
Із (~\ref{30}) та (~\ref{19}) отримуємо
$$L(M(t),Q(t))\Phi^{+}(t)X(t)=0.$$
Домноживши цю рівність справа на $X^{-1}(t)$ одержимо остаточний
результат:
\begin{equation}
\label{31} L(M(t),Q(t))Y(t)\Phi^{+}(t)=0
\end{equation}
для будь-якого $t\in R$.

Отже, з того, що підпростір $M^{n}(t)$ є інтегральним, випливає
виконання всіх умов, визначених теоремою для цього випадку, досить
лише рівність (~\ref{27}) записати як рівняння підпростору
$M^{n}(t)$ в параметричному вигляді
\begin{equation}
\label{32} y=\Phi^{+}(t)x,\,\,\,\,x\in R^{n},\,\,\,\,t\in R
\end{equation}
% \input vstavka_.tex
%%% вставка
Нехай виконується умова
\begin{equation}
\label{101}\forall t\in R:\,\,\,\,L(M(t),Q(t))\Phi^{+}(t)=0.
\end{equation}
Домноживши (~\ref{101}) справа на матрицю $\Phi(t)X(t)$, де $X(t)$
--- фундаментальна матриця розв'язків рівняння (~\ref{23}) з матрицею коефіцієнтів (~\ref{24}),
$X(0)=E$, отримаємо рівність
\begin{equation}
\label{102}\forall t\in R:\,\,\,\,L(M(t),Q(t))M(t)X(t)=0.
\end{equation}
Розглянемо тепер функцію
\begin{equation}
\label{103} r=y(t)-\Phi^{+}(t)X(t)c,
\end{equation}
де $y=y(t)$ --- розв'язок рівняння (~\ref{101}) такий, що
\begin{equation}
\label{104} y(t_{0})=\Phi^{+}(t_{0})X(t_{0})c,
\end{equation}
а $X(t),\,\,\,\,X(0)=E$ --- фундаментальна матриця розв'язків
рівняння
\begin{equation}
\label{105} \frac{dx}{dt}=P(t)x
\end{equation}
з матрицею коефіцієнтів (~\ref{24}).

Диференціюючи функцію (~\ref{103}), одержуємо
\begin{equation}
\frac{dr}{dt}=Q(t)y(t)-\left(\frac{d\Phi^{+}(t)}{dt}+\Phi^{+}(t)P(t)X(t)c\right)=Q(t)(y(t)-\Phi^{+}(t)X(t)c)+\nonumber
\end{equation}
\begin{equation}
+Q(t)\Phi^{+}(t)X(t)c-\left(\frac{d\Phi^{+}(t)}{dt}+\Phi^{+}(t)P(t)X(t)c\right)=\nonumber
\end{equation}
\begin{equation}
\label{106}
=Q(t)r-\left(\frac{d\Phi^{+}(t)}{dt}+\Phi^{+}(t)P(t)-Q(t)\Phi^{+}(t)\right)X(t)c\,\,.
\end{equation}
Доведемо, що другий доданок в (~\ref{106}) дорівнює 0. Справді,
враховуючи (~\ref{24}), маємо
\begin{equation}
\frac{d\Phi^{+}(t)}{dt}+\Phi^{+}(t)\left(\frac{d\Phi(t)}{dt}+\Phi(t)Q(t)\right)-Q(t)\Phi^{+}(t)=\frac{d\Phi^{+}(t)}{dt}+\nonumber
\end{equation}
\begin{equation}
+\Phi^{+}(t)\frac{d\Phi(t)}{dt}\Phi^{+}(t)+M(t)Q(t)\Phi^{+}(t)-Q(t)\Phi^{+}(t)=\left(\frac{d\Phi^{+}(t)}{dt}\Phi(t)+\Phi^{+}(t)\frac{d\Phi(t)}{dt}\right)\Phi^{+}(t)+\nonumber
\end{equation}
\begin{equation}
+M(t)Q(t)\Phi^{+}(t)-Q(t)M(t)\Phi^{+}(t)+Q(t)M(t)\Phi^{+}(t)-Q(t)\Phi^{+}(t)=\nonumber
\end{equation}
\begin{equation}
\label{107}
=\left(\frac{dM(t)}{dt}+M(t)Q(t)-Q(t)M(t)\right)\Phi^{+}(t)+
\end{equation}
\begin{equation}
+Q(t)\Phi^{+}(t)\Phi(t)\Phi^{+}(t)-Q(t)\Phi^{+}(t)=0.\nonumber
\end{equation}
Із врахуванням (~\ref{107}) рівність (~\ref{106}) набуває вигляду
\begin{equation}
\label{108}\frac{dr}{dt}=Q(t)r.
\end{equation}
Оскільки при $t=t_{0}$ функція (~\ref{103}) згідно з вибором
значення $y(t_{0})$ набуває нульового значення:
$$r(t_{0})=y(t_{0})-\Phi^{+}(t_{0})X(t_{0})c=0,$$
то з (~\ref{108}) випливає, що
\begin{equation}
\label{109}\forall \, t\in R:\,\,\,\,r(t)=0,
\end{equation}
а отже,
\begin{equation}
\label{110}\forall \, t\in R:\,\,\,\,y(t)=\Phi^{+}(t)X(t)c,
\end{equation}
де $c$ --- довільна точка простору $R^{n}$.

Оскільки (~\ref{110}) --- параметричний вираз рівняння підпростору
$M^{n}(t)$, то (~\ref{110}) означає, що з умови
\begin{equation}
\label{111}y(t)\in M^{n}(t)
\end{equation}
для значення $t=t_{0}$ випливає включення (\ref{111}) для довільного
$t\in R$. Це доводить, що умова (\ref{101}) є не лише необхідною,
але й достатньою для того, щоб підпростір $M^{n}(t)$ був
інваріантним многовидом рівняння (\ref{1}).

Це завершує доведення другого з тверджень теореми 1.

%%% конец вставки
Нехай підпростір $M^{m-n}(t)$ є інваріантним многовидом рівняння
(~\ref{1}). Розглянемо розв'язки $y=y(t)$ рівняння (\ref{1}),
визначені формулою
$$y(t)=Y(t)(E-M(0))c,$$
де $c$ --- довільна точка $R^{m}$.

Із того, що
$$M(0)y(0)=M(0)Y(0)(E-M(0))c=M(0)(E-M(0))c$$
випливає, що $y(0)\in M^{m-n}(0)$, а тому $\forall t\in
R:\,\,\,\,y(t)\in M^{m-n}(t)$. Це доводить, що.
\begin{equation}
\label{32'}\forall t\in R:\,\,\,\,M(t)Y(t)(E-M(0))c=0.
\end{equation}
Диференціюючи цю рівність, отримуємо
\begin{equation}
\forall t\in
R:\,\,\,\,\left(\frac{dM(t)}{dt}+M(t)Q(t)-Q(t)M(t)\right)Y(t)(E-M(0))c+\nonumber
\end{equation}
\begin{equation}
\label{33}+Q(t)M(t)Y(t)(E-M(0))c=L(M(t),Q(t))Y(t)(E-M(0))c=0.
\end{equation}
Точки $y=(E-M(0))c$ визначають підпростір $M^{m-n}(0)$, а тому
рівняння
\begin{equation}
\label{34}y=Y(t)(E-M(0))c\,\,\,\,c\in R^{m},\,\,\,\,t\in R
\end{equation}
визначає підпростір $M^{m-n}(t)$ в параметричній формі, а тому
рівність (~\ref{33}) означає, що.
\begin{equation}
\label{35}\forall \, t\in R:\,\,\,\,\ker \,L(M(t),Q(t))\supset
M^{m-n}(t)
\end{equation}
Отже, включення (~\ref{35}) є необхідною умовою для того, щоб
підпростір $M^{m-n}(t)$ був інваріантним многовидом рівняння
(~\ref{1}).

Нехай кожен із підпросторів $M^{n}(t)$ та $M^{m-n}(t)$ є
інваріантним многовидом рівняння (~\ref{1}). Тоді з рівності
(~\ref{30}) випливає, що
\begin{equation}
\label{36} \forall \, t\in R:\,\,\,\,\ker\,L(M(t),Q(t))\supset
Y(t)M(0)c,
\end{equation}
де $c$ --- довільна точка $R^{m}$. Більше того, оскільки рівняння
$$y=Y(t)(E-M(0))c\,\,\,\,c\in R^{m},\,\,\,\,t\in R$$
визначає підпростір $M^{n}(t)$ в параметричній формі, то із
(~\ref{36}) випливає включення
\begin{equation}
\label{37}\forall \,t\in R:\,\,\,\,\ker \,L(M(t),Q(t))\supset
M^{n}(t)
\end{equation}
Отже, якщо кожен із підпросторів $M^{n}(t)$ та $M^{m-n}(t)$ є
інваріантним многовидом рівняння (~\ref{1}), то
\begin{equation}
\label{38}\forall \,t\in R:\,\,\,\,\ker \,L(M(t),Q(t))\supset
(M^{n}(t)\cup M^{m-n}(t))
\end{equation}
Оскільки $\rank\,\,M^{n}(t)=n,\,\,\,\,\rank \,M^{m-n}(t)=m-n$ та
$M^{n}(t)\cap M^{m-n}(t)=\{0\}$, то об'єднання в правій частині
формули (~\ref{38}) містить в собі базу простору $R^{m}$. Тому з
(~\ref{31}) випливає, що
\begin{equation}
\label{39}\forall \,t\in \mathbb{R} :\,\,\,\,\ker
\,L(M(t),Q(t))\supset R^{m}.
\end{equation}
Оскільки $ \ker \,L(M(t),Q(t))\in M_m (\mathbb{R})$
 то включення (~\ref{39}) можливе лише тоді, коли
\begin{equation}
\label{40}\ker \,L(M(t),Q(t))=0. %(??? образ оператора)
\end{equation}
Отже, умова (~\ref{40}) є не лише достатньою, а й необхідною для
того, щоб кожен із підпросторів $M^{n}(t)$ та $M^{m-n}(t)$ був
інваріантним многовидом рівняння (\ref{1}).
\vspace*{6.0mm}  \\
%%%%%%
%------------------------------------

%------------------------------------
%\qquad \forall \, t \in \mathbb{R}

%\left( \frac{dM(t)}{dt}+M(t)Q(t)-Q(t)M(t)\right)

%\nonumber \end{equation}    ----      \rank \,

%\begin{equation} \label{№№} =========
%\end{equation}

% (\ref{1})

% розв'язки    з'єд

%\quad  \qquad  \;  \, \/

%\vspace*{6.0mm}  \\

\noindent 1. \emph{Боголюбов Н.\/} \, О некоторых статистических
методах в математической физике: Изд. АН УССР, 1945. --- Николай
Николаевич Боголюбов: Собр. научн. тр. в 12 томах, т. IV.
--- М.: Наука, 2006. --- 432~с. \\
2.  \emph{Боголюбов Н. Н., Митропольский Ю. А.\/} \, Асимптотические
методы в теории нелинейных колебаний: М., Физматгиз, 1963.
--- Николай Николаевич Боголюбов: Собр. научн. тр. в 12 томах, т. III.
--- М.: Наука, 2005. --- 605~с. \\